\documentclass[12pt]{article}
\input epsf
\usepackage{amsmath}
\usepackage{amssymb}
\usepackage{amscd}
\usepackage{epsfig}

\def\no{\noindent}
\newtheorem{lem}[equation]{Lemma}

\newtheorem{thm}[equation]{Theorem}
\newtheorem{defn}{Definition}
\newtheorem{ques}[equation]{Question}
\newtheorem{rem}[equation]{Remark}

\newcommand{\R}{\mathbb R}

\newcommand{\Z}{\mathbb Z}

\renewcommand{\H}{\mathbb H}

\newcommand{\qed}{{\unskip\nobreak\hfil
        \penalty50\hskip1em\hbox{}\nobreak\hfil
        $\square$\parfillskip=0pt\finalhyphendemerits=0 \par}}
\newcommand{\proof}{\no{\em Proof.\ }}
\def\D{\partial}
\newcommand{\al}{\alpha}

\def\Del{\Delta}
\def\eps{\epsilon}
\def\ga{\gamma}
\def\be{\beta}
\def\Ga{\Gamma}

\def\lang{\langle}
\def\<{\lang}
\def\>{\rangle}

\def\Om{\Omega}
\def\F{{\mathcal F}}
\def\C{{\mathbb C}}
\def\P{{\mathbb P}}

\def\Si{\Sigma}

\def\t{\tilde}
\def\f{\t{f}}

\def\hook{\hookrightarrow}

\begin{document}

\title{On  normal subgroups in the fundamental groups of complex surfaces}
\author{
Michael Kapovich\thanks{Supported by NSF grant
DMS-96-26633}}
\date{August 17, 1998}
\maketitle

\begin{abstract}
We show that for each aspherical compact complex surface $X$ 
whose fundamental group $\pi$ fits into a short exact sequence
$$
1\to K \to \pi  \to \pi_1(S) \to 1
$$
where $S$ is a compact hyperbolic Riemann surface 
and the group $K$ is finitely-presentable, there is a  complex structure on $S$ and a 
nonsingular holomorphic fibration $f: X\to S$ which induces the above 
short exact sequence. In particular, the fundamental groups of compact
 complex-hyperbolic surfaces cannot fit into the above short 
 exact sequence. 
As an application we give the first 
example of a non-coherent uniform lattice in $Isom(\H^2_{\C})$.

\end{abstract}

\section{Introduction}

The goal of this paper is threefold: 

(a) We will establish a restriction on the fundamental groups of compact 
aspherical complex surfaces. 

(b) We find the first examples of incoherent uniform lattices in $PU(2,1)$. 

(c) We show that the answer to the Question \ref{fibe} below is negative 
in the class of uniform lattices in $PU(2,1)$. 

\begin{ques}
\label{fibe}
Is there a Gromov-hyperbolic group $\pi$ 
which fits into a short exact sequence:
$$
1\to K \to \pi \to Q\to 1
$$ 
where $K$  and $Q$ are closed hyperbolic surface groups? 
\end{ques}

Suppose that $X$  is an  aspherical compact complex surface whose fundamental group $\pi$ 
fits into a short exact sequence
$$
1\to K \to \pi  \to Q=\pi_1(S) \to 1
$$
where $S$ is a compact hyperbolic Riemann surface 
and the group $K$ is finitely-presentable. The main theorem of this paper is 

\begin{thm}
\label{main}
Under the above assumptions there a complex structure on $S$ and a 
nonsingular holomorphic fibration $f: X\to S$ which induces the above short 
exact sequence. 
\end{thm}

\begin{rem}
Actually in Theorem \ref{main} it is enough to assume that $Q$ is a 
torsion-free group with nonzero $\be^{(2)}_1(Q)$, the 1-rst 
$L_2$-Betti number. On the other hand, in this case we have to assume that 
$X$ is K\"ahler.   Our proof also works under the assumption 
 that the group $K$ is of the type $FP_2$. 
\end{rem}

After proving Theorem \ref{main} I have learned that J.~Hillman 
\cite{Hillman(1998)} 
proved the same result under stronger assumption that $K$ is 
the fundamental group of a compact Riemann surface. Our methods 
seem to be completely different except application of 
the result of \cite{Arapura-Bressler-Ramachandran(1993)}. Later  
it turned out that the same result as Hillman's was 
independently proven by D.~Kotschick.

{\bf Acknowledgements.} I am grateful to Bill Goldman, Janos Kollar, 
Mohan Ramachandran and Domingo Toledo for discussions of issues related 
to this paper and complex-hyperbolic manifolds in general.

\section{Milnor fibration} 
\label{s1}

Let $f: \C^2 \to \C$ be a nonconstant 
holomorphic function, we assume that $0\in \C^2$ is a critical point 
of $f$. Let $S_{\eps}=S_{\eps}(0)$ 
be a sufficiently small metric sphere in $\C^2$ centered at the origin. 
Let $B_{\eps}(0)$ denote the closed $\eps$-ball centered at the origin. 
Let $K:= f^{-1}(0)\cap S_{\eps}$, this is a smooth knot (or link) in the 
3-sphere. 
The {\em Milnor fibration} $\phi : S_{\eps}-K \to S^1$ associated with $f$ is 
defined as $\phi(z,w)= f(z,w)/|f(z,w)|$, see \cite[\S 4]{Milnor(1968)}. 

Below we list  some properties of $\phi$ (see \cite[\S 4]{Milnor(1968)}, 
\cite{Eisenbud-Neumann}):  

(a) If $\eps$ is sufficiently small then $\phi$ 
determines a smooth fibration of $S_{\eps}- K$ over $S^1$. 

(b) Fibers of $\phi$ are connected provided that the germ of $f$ at zero is 
reduced, otherwise $\phi$ will have disconnected fibers. 

(b) The knot (link) $K$ is distinct from a single unknot in $S^3$ unless 
the germ of $f$  at $0$ is isomorphic  to $((z,w)\mapsto z^p, (0,0))$. 

(c) If $K$ not an unknot, then each component of $\phi^{-1}(t), t\in S^1,$ 
is not simply-connected. 

(d) Let $r>0$ be sufficiently small. Consider $s\in C_{r}(0)$, a 
point on the unit circle in $\C$ centered at zero. Let 
$\F_{\eps,s}:= f^{-1}(s)\cap B_{\eps}$. 
The two surfaces $\F_{\eps,s}$ and 
$F_{\eps,s}= \phi^{-1}(s/|s|) - f^{-1}(B_{r}(0))$ share common boundary. 
There exists an isotopy of $\F_{\eps,r}$ to $F_{\eps,r}$ within $B_{\eps}(0)$ 
which is the identity on the boundary of each surface.

\section{Multicurves}

\begin{defn}
Let $f: X\to S$ be a nonconstant proper holomorphic map from a connected 
complex surface $X$ to a Riemann surface (i.e. complex curve) $S$. We 
will say that $f$ is a {\bf nonsingular holomorphic fibration} if $f$ is a 
submersion. 
\end{defn}

Clearly the mapping $f$ as above is a real-analytic fibration, however in most 
cases it does  not determine a locally trivial holomorphic bundle. If 
$f$ is not a submersion we will still think of it as a {\bf singular} 
fibration, we shall use the notation $\F_t$ to denote the fiber $f^{-1}(t)$ 
of $f$ over $t\in S$.

\begin{defn}
Let $f: X\to D^2$ be a nonconstant proper holomorphic map with connected  
fibers where $X$ is a 2-dimensional complex surface and $D^2$ is the unit disk 
in $\C$. We assume that the origin is the only critical value of $f$.  
The singular fiber $C= f^{-1}(0)$ is called a {\bf multicurve} 
if it is a smooth curve of the multiplicity $> 1$. In other words, the germ of $f$ at each point $c\in C$ 
is equivalent to the map $(z,w) \mapsto z^n, n> 0, z, w\in \C$. The number 
$n$ is the multiplicity of $C$. 
\end{defn}

Let $t\in D^2 - 0$. Define the maps 
$$
\iota_*: H_2(f^{-1}(t)) \to H_2(X)\cong H_2(C)
$$
$$
\iota_{\#} : \pi_1(f^{-1}(t)) \to \pi_1(X) \cong \pi_1(C)
$$
induced by the inclusion $\iota: f^{-1}(t)\hook X$. 
 
\begin{lem}
\label{multicurve}
If $C$ is a multicurve then the map $\iota_*$ is not surjective. Assume that  
$C$ is a non-simply-connected multicurve.   Then the map $\iota_{\#}$ 
is not onto. 
\end{lem}
\proof Consider $Y= f^{-1}(D)\subset X$ where 
$D=\{z\in \C: |z|\le |t|\}$ is the closed 
disk in $D^2$ containing $t$. The inclusion $Y \hook X$ is a 
homotopy-equivalence so we restrict our discussion to $Y$. The 
map $C\hook Y$ is a homotopy-equivalence, thus 
the fundamental class of $C$ generates $H_2(Y)$. The dual generator 
of $H_2(Y , \D Y)$ 
is represented by 2-disk $\Del\subset Y$ which is transversal to the 
fibers of $f$ and 
$\D \Del \subset \D Y$. Since $C$ is a multicurve, the algebraic 
intersection number $[f^{-1}(t)] \cdot [\Del] =n > 1$, 
where $n$ is the multiplicity of $C$. Thus $[f^{-1}(t)]= n [C]$ 
which proves the first assertion. 

The map $\iota_{\#}$ is injective (since $\iota$ is homotopic to a 
covering $f^{-1}(t) \to C$). Thus 
$n=| \pi_1(C) : \iota_{\#}(\pi_1(f^{-1}(t)))|$, 
this proves the second assertion. \qed

\section{Proof of the main theorem}

 If $\pi_1(X)$ 
fits into short exact sequence
$$
1\to K \to \pi  \to Q=\pi_1(S) \to 1
$$
where $S$ is a hyperbolic Riemann surface 
then it follows from Kodaira's classification theorem that $X$ is a 
complex-algebraic surface. If $X$ is assumed to be K\"ahler, 
$Q$ torsion-free  and $\be^{(2)}_1(Q)\ne 0$, then $Q$ is the fundamental 
group of a hyperbolic Riemann surface, moreover if $\t{X}$ is 
the covering of $X$ corresponding to $K$ then there is a 
discrete faithful conformal action of $Q$ on $\H^2$ and a 
$Q$-equivariant proper holomorphic  map
$$
\t{f}: \t{X} \to \H^2
$$
with connected fibers (see \cite{Arapura-Bressler-Ramachandran(1993)}). 
In particular, the  projection $\pi_1(X)\to Q$ 
is induced by a holomorphic map $f: X\to S$, for the  
complex structure on $S$ given by $\H^2/Q$.   

The $i$-th $L_2$-Betti number $\be_i^{(2)}(G)$ of a finitely presentable   
group $G$ is the dimension of the $i$-th reduced $L_2$-cohomology 
group $\overline\ell_2 H^i(G)$, we refer the reader to  
\cite[Chapter 8]{Gromov(1993)} and \cite{Arapura-Bressler-Ramachandran(1993)} 
for the precise definitions. For our purposes it is enough to know that 
$\be_i^{(2)}(Q)>0$ for each 2-dimensional finitely presentable group $Q$ 
provided that $\chi(Q)< 0$  (see \cite[Chapter 8]{Gromov(1993)}). 
In particular, if $Q$ is the fundamental group of a hyperbolic 
Riemann surface of finite type then $\be^{(2)}_1(Q) >0$. Thus, 
in any case we have a holomorphic map $f: X\to S$.

We start the proof with the simple case 
when $f$ is a {\em holomorphic Morse function}, i.e. the germ of $f$ at each 
critical point is equivalent to $(z,w) \mapsto zw$. The proof in this 
case is easier and it illustrates the idea of the proof in the general case. 

Let $d$ denote the hyperbolic metric on the unit disk in $\C$.  
We will suppose that the origin $0$ is a regular value of $\f$. Direct 
computations show that the function 
$$
\ga: x\mapsto d(0, \f(x))
$$
is a real Morse function on $\t{X}$ away from $\f^{-1}(0)$ and 
the Morse index of $\ga$ at each critical point in $\t{X}- \f^{-1}(0)$ is two. 
It is clear that $r\in \R_+$ is a critical value of $\ga$ 
if and only if there is a critical value $z\in \H^2$ of $\f$ within the 
distance $r$ from the origin. Let $\F$ denote the generic fiber of $\f$. 
Thus the space $\t{X}$ is obtained by attaching 2-handles to $\F\times D^2$. 
Each singular fiber of $\f$ is obtained from $\F$ by ``pinching'' a certain 
collection of disjoint simple loops. Since $\t{X}$ is aspherical, 
each of these loops is homotopically nontrivial and no two such loops are 
homotopic to each other. (Otherwise $\t{X}$ contains a rational curve 
which then lifts to a homologically nontrivial 2-cycle in the universal 
cover of $X$.) 

We now claim that the group $\pi_1(\t{X})$ is finitely generated 
but not finitely presentable. Our proof 
follows an argument of Bestvina and Brady \cite{Bestvina-Brady}. 
Since $\t{X}$ is obtained from $\F\times D^2$ by attaching only 2-handles, 
the fundamental group of $\t{X}$ is the quotient of $\pi_1(\F)$. Recall that    
$\pi_1(\t{X})$ is finitely presentable, the epimorphism
$$
\pi_1(\F)\to \pi_1(\t{X})
$$
determines a finite generating set for $\pi_1(\t{X})$ (i.e. the generators 
of $\pi_1(\F)$).

\begin{lem}
Let $G$ be a finitely presentable group and $\{y_1,..., y_m\}$ 
be a finite generating set for $G$. Then there is a finite number 
of relators $R_1,..., R_k$ such that $\< y_1,..., y_m| R_1,..., R_k\>$ 
is a presentation of $G$. 
\end{lem}
\proof Let $\< x_1,...,x_s| Q_1,..., Q_n\>$ be a finite presentation 
of $G$. There is a finite sequence of {\em Tietze  transformations} 
(see for instance \cite[\S 1.5]{Karrass-Magnus-Solitar(1966)}) which transform 
the generating set $X=\{x_1,...,x_s\}$ to $Y=\{y_1,..., y_m\}$, 
simultaneously they transform system of relators $Q_1,..., Q_n$ for 
$X$ to a system of relators for $Y$. On each step a finite presentation 
is transformed to a finite presentation. Hence, in the end we get a 
finite system of relators $R_1,..., R_k$ for the generating set $X$. \qed

\medskip
Therefore there are finitely many elements $\al_1,...,\al_n$ of $\pi_1(\F)$ 
which 
normally generate  the kernel $Ker(\phi)$ of
$$
\phi: \pi_1(\F)\to \pi_1(\t{X})
$$
We shall identify $\al_j$ and the corresponding loops on $\F$. 
Thus there is a closed metric disk $D$ centered at the origin in $\H^2=\t{S}$ 
such that each $\al_j, j=1,...,n$, is contractible in $U=\f^{-1}(D)$. This 
implies that each $\al\in Ker(\phi)$ is contractible in $\f^{-1}(D)$. We will 
assume that the boundary of $D$ contains no critical values of $\ga$. 
However we have 
infinitely many critical values of $\f$ outside of the disk $D$. Let $z$ 
be one of them and $D'$ be a closed topological disk in $\H^2$ which 
contains both 
$D$ and $z$ and does not contain any critical  
values of $\f$ which are not in $\{z\} \cup D$. Homotopically the Morse 
surgery corresponding to $z$ amounts to attaching 2-cells along certain 
loops $\al\subset  \F$. Thus $\al\in Ker(\phi)$, which implies that 
 $\al$ is contractible 
in $U$. It follows that we get an immersed homotopically nontrivial  2-sphere  
$\zeta \subset \f^{-1}(D')$. 
The space $\t{X}$ is obtained from $\f^{-1}(D')$ 
by attaching only 2-handles, thus the homotopy 
class $[\zeta]$ is nontrivial in $\pi_2(\t{X})$ 
which contradicts asphericity of $\t{X}$. 
This concludes the proof in the case when $\f$ is a complex Morse function. 

\begin{rem}
J.~Kollar had suggested an argument which reduces the general case to the 
case of holomorphic Morse function provided that no irreducible component 
of each singular fiber of $\f$ has multiplicity $> 1$. Namely, perturb 
$\f: \t{X}\to \H^2$ in a $Q$-equivariant manner to a smooth 
map $g: \t{X}\to \H^2$ with connected fibers so that:

(a) The sets of critical values of $g$ and $\f$ are equal. 

(b) If $s$ is a critical value of $g$ (and $\f$) and $C_s\subset \H^2$ 
is a small circle around $s$ then the 3-manifolds $g^{-1}(C_s)$, 
$\f^{-1}(C_s)$ are homeomorphic.  

(c) The mapping $g$ is a holomorphic Morse function near each singular fiber. 

Then apply the same arguments as before to the function $g$ to conclude that 
neither $\f$ nor $g$ has critical points. 

However, technically it seems (at least to me) easier to apply the direct 
topological arguments below than to analyze the special case when a 
singular fiber of $\f$ has an irreducible component of multiplicity $>1$. 
\end{rem}

\medskip
We now consider the general case. We  will run essentially the same arguments 
as in the case of holomorphic Morse function. 
Let $\Si= \Si(\t f)$ denote the set of critical values of the 
holomorphic 
function $\t f$, $\t S':= \t S- \Si$ and $\t X':= \t{f}^{-1}(S')$.

\begin{lem}
\label{0}
(1) The fundamental group of a generic fiber $\F$ of $\f$ maps 
onto $K=\pi_1(\t{X})$. (2) No singular fiber 
of $\f$ is a multicurve, i.e. a singular fiber of $\f$ cannot be a smooth 
complex curve. 
\end{lem}
\proof  The restriction $\t f'$ of $\t f$ to $\t X'$ is a (nonsingular) fibration 
with connected 
fibers, thus $\pi_1(\F)$ is the kernel of the homomorphism
$$
\pi_1(\t f'): \pi_1(\t X') \to \pi_1(\t S')
$$
In particular, the subgroup $\pi_1(\F)$ is normal in $\pi_1(\t X')$. 
For each  puncture $s_i \in \Si$ choose a small loop on $\t S'$ 
going once around $s_i$ and choose a homeomorphic lift $\ga_i$  of this loop 
to $\t X'$. Then the group $\pi_1(\t X')$ is generated by 
$\pi_1(\F)$ and by the loops $\ga_i, s_i \in \Si$. Let $D_{s_i}$ 
denote a small metric disk on $\H^2$ centered at 
$s_i\in \Si$ (so that $D_{s_i} \cap \Si = \{s_i\}$). 
If for some $s_i$ the fundamental group of $\pi_1(\D \t f^{-1}(D_{s_i}))$ does 
not map onto   $\pi_1(\t f^{-1}(D_{s_i}))$ then it is true for infinitely 
many points $s\in \Si$ (all the points in the $Q$-orbit of $s_i$), thus 
the group $K$ cannot be  finitely generated. Thus the   map
$$
\pi_1(\t X') \to \pi_1(\t X)
$$
is onto. Since $\ga_i$-s belong to the kernel of this map we conclude that 
the group $\pi_1(\F)$ maps onto $\pi_1(\t X)$. 

If $\F_{s_i}= \f^{-1}(s_i), s_i\in \Si$, is a multicurve then 
$$
\pi_1(\D \t f^{-1}(D_{s_i}))\to \pi_1(\t f^{-1}(D_{s_i}))= \pi_1(\F_{s_i})
$$
is not onto (Lemma \ref{multicurve}), which contradicts our assumptions. This proves the second assertion of Lemma. \qed

Now suppose that $f: X\to S$ is not a nonsingular holomorphic fibration.  
Thus the map $\f$ has at least one fiber which is not a 
smooth complex curve. (By Lemma \ref{0} each singular fiber has to be of 
this type.) 
Our goal is to show that this  assumption leads to a 
contradiction. 
Let $T\subset \H^2$ be a locally finite embedded tree whose vertex set is 
$\Si$ (this tree of course is not $\pi_1(S)$-invariant). We can 
assume that edges of $T$ are geodesics in $\H^2$. 
For each vertex $s\in \Si$ of $T$ we choose a small closed 
metric disk $D_s$ centered at $s$ such that $D_s \cap T$ is equal to the 
intersection of $D_s$ and open edges of $T$ emanating from $s$. If 
$T'\subset T$ 
is a subtree then $N(T')$ will denote the union of $T'$ and disks $D_s$ 
for those vertices $s$ of $T$ which belong  to $T'$. Let 
$Y(T'):= \f^{-1}(N(T'))$.

Since $\f$ is a smooth fibration away from singular fibers it 
follows that the inclusion 
$$
Y(T) \hook \t{X}
$$
is a homotopy-equivalence. Therefore we restrict our attention 
to the topology of $Y(T)$. 

Let $T'$ be a finite subtree of $T$ which is the convex hull of its vertices.

\begin{lem}
\label{1}
The homomorphism 
$$
\pi_2(Y(T'))\to \pi_2(Y(T))
$$
is injective. 
\end{lem}
\proof It is enough to prove this assertion for the lifts of $Y(T')$, $Y(T)$ to 
the universal cover of $X$. Since $X$ is aspherical, its universal cover 
$\hat X$  cannot contain compact complex curves, hence the lift of 
$\f^{-1}(t), t \in T- \Si$ to $\hat X$ is a noncompact surface. Therefore this 
lift has trivial $H_2$ and the assertion follows from the 
Meyer-Vietors sequence. \qed

\medskip
Let $s\in \Si- T'$ be a vertex of $T$ which is connected to $T'$ by an edge 
$[ss']$, $s'\in \Si\cap T'$. Note that the inclusions
$$
Y(T') \hook Y(T'\cup [s's)), \quad Y(s) \hook Y([ss')) 
$$ 
are homotopy-equivalences. Here and in what follows $[ss')$ denotes the half-open edge connecting $s$ to $s'$: $s\in [ss'), s'\notin [ss')$.

\begin{lem}
\label{2}
Suppose that $\pi_1(Y(T'))\to \pi_1(Y(T))$ is a monomorphism. Then \newline  
$\pi_2(Y(T'\cup [ss']))\ne 0$. 
\end{lem}
\proof Let $t\in [ss']$ be the midpoint. Then there is a subsurface 
$\F'\subset \F_t$ such that:

(a) No boundary loop of $\F'$ is nil-homotopic in ${\F}_t$. 

(b) The image of $\pi_1(\F')$ in $\pi_1(Y([ss')))$ is trivial. 

\noindent The subsurface $\F'$ appears as follows: let $p\in \F_{s}$ be a 
singular point, then $\F'$ is a part of $\F_t$ corresponding to the Milnor 
fiber in $S_{\eps}(p)$, see section \ref{s1}. If a boundary loop of $\F'$ 
is nil-homotopic in $\F_t$ then $\t{X}$ contains a rational complex curve 
which contradicts the assumption that $\pi_2(X)=0$. 

\medskip
Therefore, the assumption of Lemma implies that the image of   
$\pi_1(\F')$ in  \newline $\pi_1(Y(T'\cap (ss']))$ is trivial. 
Consider the total lift $\widehat{\F}'$ of $\F'$ to the universal cover $\hat{X}$ of $X$, then $\widehat{\F}'$ is contained in the lift $\hat\F_t$ of $\F_t$ to $\hat{X}$. 
Note that no component of $\widehat{\F}_t- \widehat{\F}'$ is bounded 
(otherwise 
after degeneration of $\widehat{\F}_s$ to a singular fiber  we will get a 
compact complex curve in $\hat X$ which is impossible). 
If $\F'$ is not a planar surface then $\widehat{\F}'$ contains a non-separating 
loop, otherwise a component of $\D\widehat{\F}'$ is not nil-homologous 
in $\widehat{\F}_t$. 
In the both cases we apply Meyer-Vietors arguments to get 
a homologically nontrivial spherical cycle in $\hat{Y}(T'\cup [ss'])$, thus 
$\pi_2(Y(T'\cup [ss']))\ne 0$. \qed

\medskip
Since $K$ is assumed to be finitely-presentable, there are finitely many
elements $\al_i \in \pi_1(\F)$ which normally generate the kernel of
$\pi_1(\F) \to \pi_1(Y)$.  Thus there is a finite subtree $T'\subset T$ such
that all the loops $\al_i$ are nil-homotopic in $Y(T')$. Since
$\pi_1(Y(T'))$ maps onto $\pi_1(Y(T))$ (Lemma \ref{0}) it follows that 
$\pi_1(Y(T'))\to \pi_1(Y(T))$ is an isomorphism. Hence for an edge 
$[ss']$ of $T$ which has one vertex in $T'$ and the other vertex in $T- T'$ 
we have:
$$
\pi_2(Y(T'\cup [ss']))\ne 0$$
(according to Lemma \ref{2}). Now we apply Lemma \ref{1} to conclude that $\pi_2(Y(T))\ne 0$. 
However $\pi_2(Y(T))\cong \pi_2(\t X)=0$ since $X$ is aspherical. This 
contradiction proves Theorem \ref{main}. \qed

\section{Complex-hyperbolic surfaces}

Let $B\subset \C^2$ be the unit ball. We will give $B$ the {\em Kobayashi 
metric}, this metric can be described as follows. Let $p,q\in B$ be 
distinct points, there is a unique complex line 
$L\subset \C^2$ so that $p,q\in B\cap L$. Now identify $B\cap L$ with 
the hyperbolic plane $\H^2$ where the curvature is normalized to be $-1$. 
Finally let $d(p,q):= d_{\H^2}(p,q)$. Then the {\em complex-hyperbolic plane} 
$\H^2_{\C}$ is the unit ball $B$ with the Kobayashi distance $d$.  
It turns out that the Kobayashi distance $d$ is induced by a Riemannian 
metric $\rho$ on $B$. Below we list some properties 
of the complex-hyperbolic plane $\H^2_{\C}$, we refer to 
\cite{Goldman(1998)}, \cite{Hirz}, \cite{Deligne-Mostow}, \cite{Yau(1977)} for 
detailed discussion.  

\medskip
(a) $\rho$ is K\"ahler. 

(b) The sectional curvature of $\rho$ is pinched between the constants 
$-1$ and  $-1/4$. 

(c) The group of biholomorphic automorphisms of $B$ equals the 
identity component in the 
isometry group of $\H^2_{\C}$ which is isomorphic to 
the $PU(2,1)$ so that $B$ is the symmetric space for the group $PU(2,1)$: 
$B= PU(2,1)/K$ where $K\cong U(2)$ is a maximal compact subgroup  
in $PU(2,1)$.  

(d) Let $\Ga$ be a torsion-free uniform lattice in $PU(2,1)$. The quotient 
$B/\Ga$ is a compact K\"ahler surface which is actually a smooth 
complex algebraic surface. The quotient  $B/\Ga$ is called a 
{\em complex-hyperbolic surface}. 

(e) For each compact complex-hyperbolic surface we have the following 
identity between the Chern classes: $c_1^2 = 3c_2$, i.e. 
$\chi=3\tau$ where $\chi$ is the Euler characteristic and $\tau$ is 
the signature.  

(f) If $X$ is a smooth compact complex algebraic surface for which the equality 
$c_1^2 = 3c_2$ holds, then the universal cover of $X$ is biholomorphic to 
either $\H^2_{\C}$, or $\C^2$, or the complex-projective plane $\P^2_{\C}$. 

\medskip
The key fact about complex-hyperbolic surfaces which will be used in this 
paper is the following recent theorem of K.~Liu \cite{Liu(1996)}:

\begin{thm}
Let $X$ be a compact complex-hyperbolic surface. Then $X$ does not admit 
nonsingular holomorphic fibrations over complex curves. 
\end{thm}

\section{Incoherent example} 

Recall that a group $\Ga$ is called {\em coherent} if every 
finitely-generated subgroup $\Ga'\subset \Ga$ is also finitely presentable. 
Examples of coherent groups include free groups, surface groups, 
3-manifold groups (see \cite{Scott(1973a)}) 
and certain groups of cohomological dimension 2 
(see \cite{Feighn-Handel}, \cite{McCammond-Wise}). The simplest 
example of noncoherent group is ${\mathbb F}_2 \times {\mathbb F}_2$, where 
${\mathbb F}_2$ is 
the free group on two generators. (The finitely generated infinitely 
presentable subgroup in ${\mathbb F}_2\times {\mathbb F}_2$ is the kernel 
of the homomorphism 
$\phi: {\mathbb F}_2\times {\mathbb F}_2\to \Z$ where $\phi$ maps each 
free generator of 
each ${\mathbb F}_2$ to the generator of $\Z$.)  Thus there is a uniform lattice 
in the Lie group $PSL(2,\R)\times PSL(2,\R)$ which is not coherent. 
The first example of noncoherent discrete geometrically finite subgroup of 
$Isom(\H^4)$ was constructed in \cite{KP1}, \cite{Potyagailo(1990)}. 
Later on this example was generalized in \cite{Bowditch-Mess} 
to a uniform lattice in $Isom(\H^4)$.  

As an application of the main result of this paper we show that 
certain uniform lattices in $PU(2,1)= Isom(\H^2_{\C})$ are not coherent 
(these are the first known examples of incoherent discrete subgroups of 
$PU(2,1)$). The groups which we consider were known  before 
(see \cite{Livne}, \cite{Hirz}, \cite{Deligne-Mostow}) 
however their incoherence was unknown.  

\begin{lem}
\label{coh}
Suppose that $X$  is a compact complex-hyperbolic surface whose 
fundamental group $\pi$ fits into a short exact sequence
$$
1\to K \to \pi  \to Q=\pi_1(S) \to 1
$$
where $S$ is a compact hyperbolic Riemann surface 
and the group $K$ is finitely generated. 
Then $K$ is not finitely presentable. 
\end{lem}
\proof Suppose that $K$ is finitely presentable. The surface $X$ is aspherical 
since its universal cover is the complex ball. 
Then by Theorem \ref{main} 
the projection $\pi\to Q$ is induced by a nonsingular holomorphic fibration 
of the surface $X$. On the other hand, complex-hyperbolic surfaces do not admit 
such fibrations by \cite{Liu(1996)}. \qed 

\medskip
Now we describe an example the fundamental group of a complex-hyperbolic 
surface satisfying the conditions of Lemma \ref{coh} following \cite{Livne}. 
Define automorphisms $\phi, \psi$ of the free group on three generators 
$A_1, A_2, A_3$ by
$$
\phi(A_1)= A_1 A_2 A_1 ^{-1}, \phi(A_2)= A_1A_3A_1^{-1}, \phi(A_3)= A_1
$$
$$
\psi(A_1)= (A_1A_2) A_3 (A_1A_2)^{-1}, \psi(A_2)= A_1A_2A_1^{-1}, \psi(A_3)= A_1
$$
Ron Livne \cite{Livne} constructed a uniform lattice $\Ga_{d,N}$ in $PU(2,1)$ 
with the presentation
$$
\< x, y, A_1, A_2, A_3 | \ \ xA_i x^{-1}= \phi(A_i), y A_i y^{-1}= \psi(A_i)\quad  
(1\le i\le 3), $$
$$ 
x^3= y^2 = A_1 A_2 A_3, (A_1 A_2 A_3)^{2d} = A_1^2 = 
A_2^2 = A_3^2 = (yx^{-1})^N= 1\> 
$$ 
where $(N,d)\in \{ (7,7), (8,4), (9,3), (12,2)\}$. Note that the subgroup 
$K_{d}$ generated by $A_1, A_2, A_3$ in $\Ga_{d,N}$ is normal and 
finitely generated, the quotient $\Ga_{d,N}/K_{d}$ is the hyperbolic 
triangle group
$$
\Del_N:= \< x, y| x^3 = y^2 = (yx^{-1})^N= 1\>
$$
since $N\ge 7$. Now fix a pair $(N,d)$ from the above list and let 
$\Ga:= \Ga_{d,N}, \Del:= \Del_{N}, K:= K_{d}$. 
Let $\Del'< \Del$ be a torsion-free subgroup of finite index 
and $\Ga'< \Ga$ be the pull-back of $\Del'$ to $\Ga$. Then $\Ga'$ fits 
into short exact sequence
$$
1\to K \to \Ga' \to \Del' \to 1
$$ 
The group $\Ga'$ still has torsion, so let $\pi$ be a torsion-free subgroup 
of $\Ga'$, $K':= \pi\cap K, Q:= \pi/K'$. Clearly $K'$ is finitely generated 
and $Q$ is the fundamental group of a compact hyperbolic Riemann surface. 
The group $\pi$ acts freely discretely cocompactly on $\H^2_{\C}$ and hence 
is the fundamental group  of the compact complex-hyperbolic surface 
$X= \H^2_{\C}/\pi$.  By Lemma \ref{coh} 
the group $K'$ is not finitely presentable.

\begin{rem}
Bill Goldman had told me long ago about Livne's example as a candidate for 
non-coherence, however until recently I did not know how to prove that the 
group $K$ is not finitely presentable. 
\end{rem}

Note that the group $K$ is  not geometrically finite and its limit 
set is the whole sphere at infinity of $\H^2_{\C}$ (since $K$ is 
normal in $\Ga$). 

\begin{ques}
Let $\Ga\subset  PU(2,1)$ be a finitely generated discrete subgroup whose 
limit set is not the whole sphere at infinity of $\H^2_{\C}$. Is $\Ga$ 
finitely-presentable? Is $\Ga$ geometrically finite? 
\end{ques}

\begin{rem}
There are several reasons why it is difficult to construct finitely generated 
geometrically infinite subgroups of $PU(2,1)$. One of them is the following 
result due to M.~Ramachandran:

Let $\Ga$ be a discrete subgroup of $PU(2,1)$ which does not contain 
parabolic elements and which acts cocompactly on a component $\Om_0$ 
of the domain of discontinuity $\Om(\Ga)\subset \D_{\infty}\H^2_{\C}$. Then 
 $\Ga$ is geometrically finite and $\Om_0= \Om(\Ga)$. 
 
(Instead of assuming that $\Ga$ contains no parabolic elements it is enough 
to assume that each maximal parabolic subgroup of $\Ga$ is isomorphic to 
a lattice in the 3-dimensional Heisenberg group. )
\end{rem}

\begin{ques}
Is there a compact {\bf real-hyperbolic} $4$-manifold $X$ whose 
fundamental group fits into a short exact sequence:
$$
1\to K \to \pi_1(X) \to Q\to 1
$$ 
where $K$ is finitely presentable or even a surface group and $Q$ is a  
hyperbolic surface group ?
\end{ques}

More generally:

\begin{ques}
Is there a Gromov-hyperbolic group $\pi$ 
which fits into a short exact sequence:
$$
1\to K \to \pi \to Q\to 1
$$ 
where $K$  and $Q$ are closed hyperbolic surface groups? 
\end{ques}

Note that Lee Mosher \cite{Mosher} constructed similar example when $K$ 
is a closed hyperbolic surface group and $Q$ is a free nonabelian group. 

\begin{ques}
Let $\Ga_g$ be the mapping class group of a compact surface of genus $g$. 
Is there $g$ and a finitely generated  non-free subgroup $Q$ 
of $\Ga_g$ which consists only of the identity and pseudo-Anosov elements? 
\end{ques}

Mosher's example comes from a ``Schottky-type'' subgroup $Q$ in $\Ga_g$ where 
$K$ is the fundamental group of a genus $g$ surface.

\bibliography{/u/ma/kapovich/lit}
\bibliographystyle{siam}

\noindent Michael Kapovich: Department of Mathematics, University of 
Utah,  Salt Lake City, UT 84112, USA ; kapovich$@$math.utah.edu

\end{document}